\documentclass[12pt]{article}
\usepackage{amsmath}
\usepackage{amssymb}
\usepackage{amsthm}
\usepackage{amsfonts}
\usepackage{amscd}
\usepackage[mathscr]{eucal}
\newcommand{\Z} {{\mathbb  Z}}
\newcommand{\Q}{{\mathbb  Q}}

\newcommand{\R} {{\mathbb R}}

\begin{document}
\parindent  25pt
\baselineskip  10mm
\textwidth  15cm    \textheight  23cm \evensidemargin -0.06cm
\oddsidemargin -0.01cm

\title{ { On several families of elliptic curves with arbitrary
large Selmer groups }}
\author{\mbox{}  {Fei Li and  Derong Qiu }
\thanks{ \quad E-mail:
derong@mail.cnu.edu.cn } \\
(School of Mathematical Sciences, \\
Institute of Mathematics and Interdisciplinary Science, \\
Capital Normal University, Beijing 100048, P.R.China )   }

\date{}
\maketitle
\parindent  24pt
\baselineskip  10mm
\parskip  0pt

\par   \vskip 0.6cm

\hspace{-0.6cm}{\bf 1 \ \ Introduction and Main Results}

      \par \vskip 0.4 cm

In this paper, we consider the elliptic curve $$ E =
E_{\varepsilon }: y^{2} = x ( x + \varepsilon p D ) ( x +
\varepsilon q D ), \quad ( \varepsilon = \pm 1 ) \qquad \qquad
\quad (1.1) $$ where $ p $ and  $ q $ are odd prime numbers with $
q - p = 2, $ and $ D = D_{1} \cdots D_{n} $ is a square-free
integer with distinct primes $ D_{1}, \cdots , D_{n}. $  Moreover,
$ 2 \nmid D, \ p \nmid D $ and $ q \nmid D. $ For each $ D_{i}, $
denote $ \widehat{D_{i}} = D / D_{i} \ ( \widehat{D_{1}}= 1 \
\text{if } \ D = D_{1} ). $ We write $ E = E_{+} $ if $
\varepsilon = 1 $ and $ E = E_{-} $ if $ \varepsilon = -1. $
\par \vskip 0.2 cm

Take the elliptic curve $$ E^{\prime } = E_{\varepsilon }^{\prime
}: y^{2} = x^{3} - 2 \varepsilon ( p + q ) D x^{2} + 4 D^{2} x.
\qquad \qquad \quad (1.2) $$ There is an isogeny $ \varphi $ of
degree $ 2 $ between $ E $ and $ E^{\prime} $ as follows:
$$ \varphi : E \longrightarrow E^{\prime}, \quad ( x, y ) \longmapsto
( y^{2} / x^{2}, \ y( pqD^{2} - x^{2} ) / x^{2}). $$ The kernel is $
E [ \varphi ] = \{ O, \ (0, 0) \}, $ and the dual isogeny of $
\varphi $ is $$ \widehat{ \varphi } : E^{\prime } \longrightarrow E,
\quad ( x, y ) \longmapsto ( y^{2} /4 x^{2}, \ y( 4 D^{2} - x^{2} )
/8 x^{2}) $$ with kernel $ E^{\prime} [ \widehat{ \varphi } ] = \{
O, \ (0, 0) \} $ (see [S, p.74]).
\par \vskip 0.2 cm

In this paper, we calculate the $ \varphi ( \widehat{\varphi
})-$Selmer groups $ S^{(\varphi )} (E / \Q ) $ and $
S^{(\widehat{\varphi })} (E^{\prime} / \Q ) $ of the above elliptic
curves via descent theory (see [S, Chapter X]), in particular, we
obtain that the Selmer groups of several families of such elliptic
curves can be arbitrary large. There are many literature studying
Selmer groups of elliptic curves ( see e.g., [KS], [QZ], [Sch],
[ST]). The main results in the present paper are stated as follows:
\par \vskip 0.3 cm

\hspace{-0.6cm}{\bf Part 1  \ Results for the case \ $ \varepsilon
= + 1 $ }
\par \vskip 0.3 cm

 {\bf Theorem 1.1. } \quad Let $ E = E_{+} $ be
the elliptic curve in (1.1) with $ \varepsilon = + 1. $ For each $
i \in \{ 1, \cdots , n \}, $ denote $$ \Pi _{i}^{+}(D) = \left( 1
- \left(\frac{-1}{D_{i}}\right) \right) + \left( 1 -
\left(\frac{p}{D_{i}}\right) \right) + \left( 1 -
\left(\frac{q}{D_{i}}\right) \right) + \sum _{j=1, \ j\neq i}
^{n}\left( 1 - \left(\frac{D_{j}}{D_{i}}\right) \right), $$ where
$ \left(- \right) $ is the ( Legendre ) quadratic residue symbol.
And define a function $ \rho ^{+}( D ) $ by $$ \rho ^{+}( D ) =
\sum _{i = 1} ^{n} \left[ \frac{1}{1 + \Pi _{i}^{+}(D)}\right],
$$ where [x] is the greatest integer $ \leq x. $ Then we have
\par \vskip 0.1 cm

(1) \quad There exists a subset $ T \subset \{ D_{1}, \cdots ,
D_{n} \} $ with cardinal $ \sharp T = \rho ^{+}( D ) $ such that $
S^{(\varphi )} (E / \Q ) \supset \ < T \text{mod} (\Q^{\star^{2}})
> \ \cong \left( \Z / 2 \Z \right)^{\rho ^{+}( D )}. $ In particular,
$ \text{dim}_{2} S^{(\varphi )} (E / \Q ) \geq \rho ^{+}( D ). $
\par \vskip 0.1 cm

(2) \quad If $ p \equiv 7 (\text{mod} 8 ) $ and $ D_{i} \equiv 1,
\ 7 (\text{mod} 8 ) \ ( 1 \leq i \leq n), $ then $ \text{dim}_{2}
S^{(\varphi )} (E / \Q ) \geq \rho ^{+}( D ) + 1 $ and $
S^{(\varphi )} (E / \Q ) \supset \ < T \cup \{ 2 \} \text{mod}
(\Q^{\star^{2}})
> \ \cong \left( \Z / 2 \Z \right)^{\rho ^{+}( D ) + 1 } $ for some
subset $ T \subset \{ D_{1}, \cdots , D_{n} \} $ with  $ \sharp T
= \rho ^{+}( D ). $
\par \vskip 0.2 cm

{\bf Corollary 1.2.} \quad For elliptic curve $ E = E_{+} $ as in
Theorem 1.1, we have $ S^{(\varphi )} (E / \Q ) \ \subset \ \left(
\Z / 2 \Z \right)^{n + 1}. $ Furthermore,
\par \vskip 0.1 cm

(A) \ If the following conditions hold: \\
(1) \ $ D_{i} \equiv 1 (\text{mod} 4 ), \quad 1 \leq i \leq n; $
\quad \quad (2) \ $ \left(\frac{D_{i}}{p}\right) =
\left(\frac{D_{i}}{q}\right) = 1, \quad 1 \leq i \leq n; $ \\
(3) \ $ \left(\frac{D_{j}}{D_{i}}\right) = 1, \quad 1 \leq i, \ j
\leq n $ \ and \ $ i \neq j. $ \\
Then \ $ S^{(\varphi )} (E / \Q ) \ \supset \ \left( \Z / 2 \Z
\right)^{n}, $ \ so \ $ 2^{n} \leq \sharp S^{(\varphi )} (E / \Q )
\leq 2^{n + 1 }. $
\par \vskip 0.1 cm
(B) \ If the following conditions hold: \\
(1) \ $ D_{i} \equiv 1 (\text{mod} 4 ), \quad 1 \leq i \leq n; $ \\
(2) \ $ D_{k} \equiv 5 (\text{mod} 8 ) $ for at least one $ k \in \{
1, \cdots , n \}; $ \\
(3) \ $ \left(\frac{D_{i}}{p}\right) = \left(\frac{D_{i}}{q}\right)
= 1, \quad 1 \leq i \leq n; $ \\
(4) \ $ \left(\frac{D_{j}}{D_{i}}\right) = 1, \quad 1 \leq
i, \ j \leq n $ \ and \ $ i \neq j. $ \\
Then \ $ S^{(\varphi )} (E / \Q ) \cong \ \left( \Z / 2 \Z
\right)^{n }, $ \ and \ $ \sharp S^{(\varphi )} (E / \Q ) = 2^{n}. $
\par \vskip 0.1 cm
(C) \ If the following conditions hold: \\
(1) \ $ p \equiv 7 (\text{mod} 8 ); $ \quad \quad (2) \ $ D_{i}
\equiv 1 (\text{mod} 8 ), \quad 1 \leq i \leq n; $ \\
(3) \ $ \left(\frac{D_{i}}{p}\right) =
\left(\frac{D_{i}}{q}\right) = 1, \quad 1 \leq i \leq n; $ \\
(4) \ $ \left(\frac{D_{j}}{D_{i}}\right) = 1, \quad 1 \leq
i, \ j \leq n $ \ and \ $ i \neq j. $ \\
Then \ $ S^{(\varphi )} (E / \Q ) \cong \ \left( \Z / 2 \Z
\right)^{n + 1}, $ \ and \ $ \sharp S^{(\varphi )} (E / \Q ) =
2^{n + 1}. $

\par \vskip 0.2 cm

{\bf Theorem 1.3. }\quad Let $ E^{\prime} = E^{\prime}_{+} $ be
the elliptic curve in (1.2) with $ \varepsilon = + 1. $ For each $
i \in \{ 1, \cdots , n \}, $ denote $$ \Pi _{i}^{+}(D)^{\prime } =
\left( 1 - \left(\frac{ - p \widehat{D_{i}}}{D_{i}}\right) \right)
\left( 1 - \left(\frac{- q \widehat{D_{i}} }{D_{i}}\right) \right)
+ \sum_{j = 1, \ j \neq i}^{n} \left( 1 -
\left(\frac{D_{i}}{D_{j}}\right) \right) \left( 1 - \left(\frac{p
q D_{i} }{D_{j}}\right) \right), $$ where $ \left(- \right) $ is
the ( Legendre ) quadratic residue symbol. Take a subset $ I $ of
$ Z(n) = \{ 1, \cdots , n \} $ as follows: \ $ I = \{i \in Z(n): \
D_{i} \equiv 1 (\text{mod} 8 ) \} \cup \{i \in Z(n) : \ (1 + p
\widehat{D_{i}} )(1 + q \widehat{D_{i}} ) \equiv 0 (\text{mod} 16
) \} \cup \{i \in Z(n) : \ D_{i} \equiv 3 (\text{mod} 8 ) \
\text{and} \ p \equiv 1 (\text{mod} 4 ) \} \cup \{i \in Z(n) : \
D_{i} \equiv 7 (\text{mod} 8 ) \ \text{and} \ p \equiv 3
(\text{mod} 4 ) \}, $ and define a function $ \rho ^{+}( D
)^{\prime } $ by
$$ \rho ^{+}( D )^{\prime } = \sum _{i \in I } \left[ \frac{1}{1 +
\Pi _{i}^{+}(D)^{\prime }} \right],
$$ where [x] is the greatest integer $ \leq x. $ Then there exists
a subset $ T \subset \{ D_{1}, \cdots , D_{n} \} $ with cardinal $
\sharp T = \rho ^{+}( D )^{\prime} $ such that $
S^{(\widehat{\varphi } )} (E^{\prime } / \Q ) \supset \ < T
\text{mod} (\Q^{\star^{2}})
> \ \cong \left( \Z / 2 \Z \right)^{\rho ^{+}( D )^{\prime }}. $
In particular, $ \text{dim}_{2} S^{(\widehat{\varphi } )}
(E^{\prime } / \Q ) \geq \rho ^{+}( D )^{\prime }. $
\par \vskip 0.2 cm

{\bf Corollary 1.4.} \quad For elliptic curve $ E^{\prime} =
E^{\prime}_{+} $ as in Theorem 1.3, we have $ S^{(\widehat{\varphi
})} (E^{\prime} / \Q ) \ \subset \ \left( \Z / 2 \Z \right)^{n + 3}.
$ Furthermore, if for each $ i \in \{ 1, \cdots , n  \}, \ \Pi
_{i}^{+}(D)^{\prime } = 0 \ $ and one of the following conditions hold: \\
(1) \ $ D_{i} \equiv 1 (\text{mod} 8 ); $ \quad (2) \ $ (1 + p
\widehat{D_{i}} )(1 + q \widehat{D_{i}} )  \equiv 0 (\text{mod} 16 ); $ \\
(3) \ $ \left \{
   \begin{array}{l}
  D_{i} \equiv 3 (\text{mod} 8 ), \\
  p \equiv 1 (\text{mod} 4 ); \\
\end{array}
  \right.   \qquad  \quad (4) \ \left \{ \begin{array}{l}
  D_{i} \equiv 7 (\text{mod} 8 ), \\
  p \equiv 3 (\text{mod} 4 ).
  \end{array}
  \right. $ \\
Then \ $ S^{(\widehat{\varphi })} (E^{\prime} / \Q ) \supset \
\left( \Z / 2 \Z \right)^{n + 2 }, $ \ so \ $ 2^{n + 2} \leq \sharp
S^{(\widehat{\varphi })} (E^{\prime} / \Q ) \leq 2^{n + 3}. $ \\
Furthermore, if $ \left \{
   \begin{array}{l}
p \equiv 3 (\text{mod} 4 ), \\
\alpha ( -pq )^{\prime} = 0
\end{array}
 \right.   \ \text{or} \ \left \{ \begin{array}{l}
  D - p \equiv 0, 2 (\text{mod} 8 ), \\
  \alpha ( -pq )^{\prime} = 0,
 \end{array}
 \right.  $ \\
then $ S^{(\widehat{\varphi })} (E^{\prime} / \Q ) = \left( \Z / 2
\Z \right)^{n + 3 }, $ so \ $ \sharp S^{(\widehat{\varphi })}
(E^{\prime} / \Q ) = 2^{n + 3}, $ \\
where \ $ \alpha ( -pq )^{\prime} = \sum_{i = 1}^{n} \left( 1 -
\left(\frac{-1}{D_{i}}\right) \right)
\left( 1 - \left(\frac{ - p q }{D_{i}}\right) \right). $ \\
For example, if the following conditions hold: \\
(1) \ $ D_{i} \equiv 1 (\text{mod} 8 ), \quad 1 \leq i \leq n; $ \
\quad (2) \ $ \left(\frac{D_{i}}{p}\right) +
\left(\frac{D_{i}}{q}\right) = 0, \quad 1 \leq i \leq n; $ \\
(3) \ $ p \equiv 3 (\text{mod} 4 ). $ \\
Then $ S^{(\widehat{\varphi })} (E^{\prime} / \Q ) = \left( \Z / 2
\Z \right)^{n + 3 }, $ so \ $ \sharp S^{(\widehat{\varphi })}
(E^{\prime} / \Q ) = 2^{n + 3}. $

\par \vskip 0.2 cm

{\bf Theorem 1.5. } \quad For elliptic curves $ E = E_{+} $ and  $
E^{\prime } = E^{\prime }_{+} $in (1.1) and (1.2) with $
\varepsilon = + 1, $ we have \\
(A) \ If the following conditions hold: \\
(1) \ $ D_{i} \equiv 1 (\text{mod} 4 ), \quad 1 \leq i \leq n; $ \\
(2) \ $ D_{k} \equiv 5 (\text{mod} 8 ) $ for at least one
$ k \in \{1, \cdots , n \}; $ \\
(3) \ $ \left(\frac{D_{i}}{p}\right) =
\left(\frac{D_{i}}{q}\right) = 1, \quad 1 \leq i \leq n; $ \\
(4) \ $ \left(\frac{D_{j}}{D_{i}}\right) = 1, \quad 1 \leq
i, \ j \leq n $ \ and \ $ i \neq j; $ \\
(5) \ $ p - D \equiv 0, 2 (\text{mod} 8 ). $ \\
Then \ $ S^{(\varphi )} (E / \Q ) \cong \ \left( \Z / 2 \Z
\right)^{n},  \ \left( \Z / 2 \Z \right)^{n + 2} \subset
S^{(\widehat{\varphi })} (E^{\prime} / \Q ) \subset \left( \Z / 2 \Z
\right)^{n + 3 }. $ \\
Furthermore, if $ p \equiv 3 (\text{mod} 4 )
$ or $ p - D \equiv 0 (\text{mod} 8 ), $
then  \\
$ S^{(\varphi )} (E / \Q ) \cong \ \left( \Z / 2 \Z \right)^{n},
\quad S^{(\widehat{\varphi })} (E^{\prime} / \Q ) \cong
\left( \Z / 2 \Z \right)^{n + 3 },  $  and \\
$ \text{rank}( E(\Q )) + \text{dim}_{2} \text{TS}(E/\Q)[\varphi ] +
\text{dim}_{2} \text{TS}(E^{\prime }/\Q )[\widehat{\varphi } ] = 2n
+1, $ \\
where $ \text{TS}(E/\Q) $ and $ \text{TS}(E^{\prime }/\Q )
$ are the Tate-Shafarevich groups of $ E $ and $ E^{\prime } $ over
$ \Q $ respectively (See [S, Chapter X]). \\
(B) \ If the following conditions hold: \\
(1) \ $ p \equiv 7 (\text{mod} 8 ); $ \quad \quad (2) \ $ D_{i}
\equiv 1 (\text{mod} 8 ), \quad 1 \leq i \leq n; $ \\
(3) \ $ \left(\frac{D_{i}}{p}\right) =
\left(\frac{D_{i}}{q}\right) = 1, \quad 1 \leq i \leq n; $ \\
(4) \ $ \left(\frac{D_{j}}{D_{i}}\right) = 1, \quad 1 \leq
i, \ j \leq n $ \ and \ $ i \neq j. $ \\
Then \ $ S^{(\varphi )} (E / \Q ) \cong \ \left( \Z / 2 \Z
\right)^{n + 1}, $ \quad  $ S^{(\widehat{\varphi })} (E^{\prime}
/ \Q ) \cong \left( \Z / 2 \Z \right)^{n + 3 }, $ \ and \\
$ \text{rank}( E(\Q)) + \text{dim}_{2} \text{TS}(E/\Q)[\varphi ] +
\text{dim}_{2} \text{TS}(E^{\prime }/\Q)[\widehat{\varphi } ] = 2n
+ 2. $

\par \vskip 0.4 cm

\hspace{-0.6cm}{\bf Part 2 \ Results for the case \ $ \varepsilon
= - 1 $ }
\par \vskip 0.3 cm

{\bf Theorem 1.6. } \quad Let $ E = E_{-} $ be the elliptic curve
in (1.1) with $ \varepsilon = - 1. $ For each $ i \in \{ 1, \cdots
, n \}, $ denote $$ \Pi _{i}^{-}(D) =  \left( 1 -
\left(\frac{p}{D_{i}}\right) \right) + \left( 1 -
\left(\frac{q}{D_{i}}\right) \right) + \sum _{j=1, \ j\neq i}
^{n}\left( 1 - \left(\frac{D_{j}}{D_{i}}\right) \right), $$ and
define a function $ \rho ^{-}( D ) $ by $$ \rho ^{-}( D ) = \sum
_{i = 1} ^{n} \left[ \frac{1}{1 + \Pi _{i}^{-}(D)}\right], $$
where [x] is the greatest integer $ \leq x. $ Then we have
\par \vskip 0.1 cm

(1) \quad There exists a subset $ T \subset \{ D_{1}, \cdots ,
D_{n} \} $ with cardinal $ \sharp T = \rho ^{-}( D ) $ such that $
S^{(\varphi )} (E / \Q ) \supset \ < \{ D_{i}^{\star } : \ D_{i}
\in T \}\text{mod} (\Q^{\star^{2}})
> \ \cong \left( \Z / 2 \Z \right)^{\rho ^{-}( D )}, $ \ where
$ D_{i}^{\star } \in \{ D_{i}, \ - D_{i} \}. $ In particular, $
\text{dim}_{2} S^{(\varphi )} (E / \Q ) \geq \rho ^{-}( D ). $
\par \vskip 0.1 cm

(2) \quad Furthermore, if one of the following conditions hold: \\
(a) \ $ p \equiv 7 (\text{mod} 8 ) $ and $ D_{i} \equiv 1, \ 7
(\text{mod} 8 ) \ ( 1 \leq i \leq n); $ \\
(b) \ $ p \equiv 1 (\text{mod} 8 ) $ and $ D_{i} \equiv 1, \ 3
(\text{mod} 8 ) \ ( 1 \leq i \leq n). $ \\
Then for $ T $ as in (1), we have $ S^{(\varphi )} (E / \Q )
\supset \ < 2^{\star }> \times < \{ D_{i}^{\star } : \ D_{i} \in T
\} \cong \left( \Z / 2 \Z \right)^{\rho ^{-}( D ) + 1 }, $ where $
2^{\star } \in \{ 2, \ - 2 \}. $ In particular, $ \text{dim}_{2}
S^{(\varphi )} (E / \Q ) \geq \rho ^{-}( D ) + 1. $
\par \vskip 0.2 cm

{\bf Corollary 1.7.} \quad For elliptic curve $ E = E_{-} $ as in
Theorem 1.6, we have $ S^{(\varphi )} (E / \Q ) \ \subset \ \left(
\Z / 2 \Z \right)^{n + 1}. $ Furthermore,
\par \vskip 0.1 cm

(A) \ If the following conditions hold: \\
(1) \ $ \left(\frac{p}{D_{i}}\right) =
\left(\frac{q}{D_{i}}\right) = 1, \quad 1 \leq i \leq n; $ \\
(2) \ $ \left(\frac{D_{j}}{D_{i}}\right) = 1, \quad 1 \leq i, \ j
\leq n $ \ and \ $ i \neq j. $ \\
Then \ $ S^{(\varphi )} (E / \Q ) \supset \ \left( \Z / 2 \Z
\right)^{n}, $ \ so \ $ 2^{n} \leq \sharp S^{(\varphi )} (E / \Q )
\leq 2^{n + 1 }. $
\par \vskip 0.1 cm

(B) \ If the following conditions hold: \\
(1) \ $ \left(\frac{p}{D_{i}}\right) =
\left(\frac{q}{D_{i}}\right) = 1, \quad 1 \leq i \leq n; $ \\
(2) \ $ \left(\frac{D_{j}}{D_{i}}\right) = 1, \quad 1 \leq i, \ j
\leq n $ \ and \ $ i \neq j; $ \\
(3) \ $ \left \{
   \begin{array}{l}
  p \equiv 7 (\text{mod} 8 ) \\
  D_{i} \equiv 1, 7 (\text{mod} 8 ) \ ( 1 \leq i \leq n ),
  \end{array}
  \right.   \ \text{or} \ \left \{ \begin{array}{l}
  p \equiv 1 (\text{mod} 8 ) \\
  D_{i} \equiv 1, 3 (\text{mod} 8 ) \ ( 1 \leq i \leq n ).
  \end{array}
  \right. $  \\
Then \ $ S^{(\varphi )} (E / \Q ) \cong \ \left( \Z / 2 \Z
\right)^{n + 1}, $ \ and \ $ \sharp S^{(\varphi )} (E / \Q ) =
2^{n + 1}. $
\par \vskip 0.2 cm

{\bf Theorem 1.8.}\quad Let $ E^{\prime} = E^{\prime}_{-} $ be the
elliptic curve in (1.2) with $ \varepsilon = - 1. $ For each $ i
\in \{ 1, \cdots , n \}, $ denote $$ \Pi _{i}^{-}(D)^{\prime } =
\left( 1 - \left(\frac{ p \widehat{D_{i}}}{D_{i}}\right) \right)
\left( 1 - \left(\frac{ q \widehat{D_{i}} }{D_{i}}\right) \right)
+ \sum_{j = 1, \ j \neq i}^{n} \left( 1 -
\left(\frac{D_{i}}{D_{j}}\right) \right) \left( 1 - \left(\frac{p
q D_{i} }{D_{j}}\right) \right), $$ where $ \left(- \right) $ is
the ( Legendre ) quadratic residue symbol. Take a subset $ I $ of
$ Z(n) = \{ 1, \cdots , n \} $ as follows: \ $ I = \{i \in Z(n): \
D_{i} \equiv 1 (\text{mod} 8 ) \} \cup \{i \in Z(n) : \ (1 - p
\widehat{D_{i}} )(1 - q \widehat{D_{i}} ) \equiv 0 (\text{mod} 16
) \} \cup \{i \in Z(n) : \ D_{i} \equiv 3 (\text{mod} 8 ) \
\text{and} \ p \equiv 1 (\text{mod} 4 ) \} \cup \{i \in Z(n) : \
D_{i} \equiv 7 (\text{mod} 8 ) \ \text{and} \ p \equiv 3
(\text{mod} 4 ) \}, $ and define a function $ \rho ^{-}( D
)^{\prime } $ by
$$ \rho ^{-}( D )^{\prime } = \sum _{i \in I } \left[ \frac{1}{1 +
\Pi _{i}^{-}(D)^{\prime }} \right],
$$ where [x] is the greatest integer $ \leq x. $ Then there exists
a subset $ T \subset \{ D_{1}, \cdots , D_{n} \} $ with cardinal $
\sharp T = \rho ^{-}( D )^{\prime} $ such that $
S^{(\widehat{\varphi } )} (E^{\prime } / \Q ) \supset \ < T
\text{mod} (\Q^{\star^{2}})
> \ \cong \left( \Z / 2 \Z \right)^{\rho ^{-}( D )^{\prime }}. $
In particular, $ \text{dim}_{2} S^{(\widehat{\varphi } )}
(E^{\prime } / \Q ) \geq \rho ^{-}( D )^{\prime }. $
\par \vskip 0.2 cm

{\bf Corollary 1.9.} \quad For elliptic curve $ E^{\prime} =
E^{\prime}_{-} $ as in Theorem 1.8, we have $ S^{(\widehat{\varphi
})} (E^{\prime} / \Q ) \ \subset \ \left( \Z / 2 \Z \right)^{n + 2}.
$ Furthermore, if for each $ i \in \{ 1, \cdots , n  \}, \ \Pi
_{i}^{-}(D)^{\prime } = 0 \ $ and one of the following conditions hold: \\
(1) \ $ D_{i} \equiv 1 (\text{mod} 8 ); $ \quad (2) \ $ (1 - p
\widehat{D_{i}} )(1 - q \widehat{D_{i}} )  \equiv 0 (\text{mod} 16 ); $ \\
(3) \ $ \left \{
   \begin{array}{l}
  D_{i} \equiv 3 (\text{mod} 8 ), \\
  p \equiv 1 (\text{mod} 4 ); \\
\end{array}
  \right.   \qquad  \quad (4) \ \left \{ \begin{array}{l}
  D_{i} \equiv 7 (\text{mod} 8 ), \\
  p \equiv 3 (\text{mod} 4 ).
  \end{array}
  \right. $ \\
Then \ $ S^{(\widehat{\varphi })} (E^{\prime} / \Q ) \cong \ \left(
\Z / 2 \Z \right)^{n + 2 }, $ \ so \ $ \sharp
S^{(\widehat{\varphi })} (E^{\prime} / \Q ) = 2^{n + 2}. $ \\
For example, if the following conditions hold: \\
(1) \ $ D_{i} \equiv 1 (\text{mod} 8 ), \quad 1 \leq i \leq n; $ \
\quad (2) \ $ \left(\frac{D_{i}}{p}\right) +
\left(\frac{D_{i}}{q}\right) = 0, \quad 1 \leq i \leq n. $ \\
Then $ S^{(\widehat{\varphi })} (E^{\prime} / \Q ) = \left( \Z / 2
\Z \right)^{n + 2 }, $ so \ $ \sharp S^{(\widehat{\varphi })}
(E^{\prime} / \Q ) = 2^{n + 2}. $
\par \vskip 0.2 cm

{\bf Theorem 1.10 } \quad For elliptic curves $ E = E_{-} $ and  $
E^{\prime } = E^{\prime }_{-} $in (1.1) and (1.2) with $ \varepsilon
= -1, $ we have \\
(A) \ If the following conditions hold: \\
(1) \ $ D_{i} \equiv 1 (\text{mod} 4 ), \quad 1 \leq i \leq n; $ \\
(2) \ $ D_{k} \equiv 5 (\text{mod} 8 ) $ for at least one
$ k \in \{1, \cdots , n \}; $ \\
(3) \ $ \left(\frac{D_{i}}{p}\right) =
\left(\frac{D_{i}}{q}\right) = 1, \quad 1 \leq i \leq n; $ \\
(4) \ $ \left(\frac{D_{j}}{D_{i}}\right) = 1, \quad 1 \leq
i, \ j \leq n $ \ and \ $ i \neq j; $ \\
(5) \ $ p - D \equiv 2, 4 (\text{mod} 8 ). $ \\
Then \ $ S^{(\varphi )} (E / \Q ) \cong \ \left( \Z / 2 \Z
\right)^{n}, \quad S^{(\widehat{\varphi })} (E^{\prime} / \Q ) \cong
\left( \Z / 2 \Z \right)^{n + 2},  $  and \\
$ \text{rank}( E(\Q)) + \text{dim}_{2} \text{TS}(E/\Q)[\varphi ] +
\text{dim}_{2} \text{TS}(E^{\prime }/\Q)[\widehat{\varphi } ] = 2n. $ \\
(B) \ If the following conditions hold: \\
(1) \ $ p \equiv 1, 7 (\text{mod} 8 ); $ \quad \quad (2) \ $ D_{i}
\equiv 1 (\text{mod} 8 ), \quad 1 \leq i \leq n; $ \\
(3) \ $ \left(\frac{p}{D_{i}}\right) =
\left(\frac{q}{D_{i}}\right) = 1, \quad 1 \leq i \leq n; $ \\
(4) \ $ \left(\frac{D_{j}}{D_{i}}\right) = 1, \quad 1 \leq
i, \ j \leq n $ \ and \ $ i \neq j. $ \\
Then \ $ S^{(\varphi )} (E / \Q ) \cong \ \left( \Z / 2 \Z
\right)^{n + 1}, $ \quad  $ S^{(\widehat{\varphi })} (E^{\prime}
/ \Q ) \cong \left( \Z / 2 \Z \right)^{n + 2 }, $ \ and \\
$ \text{rank}( E(\Q)) + \text{dim}_{2} \text{TS}(E/\Q)[\varphi ] +
\text{dim}_{2} \text{TS}(E^{\prime }/\Q)[\widehat{\varphi } ] = 2n
+ 1. $
\par \vskip 0.2 cm

{\bf Remark 1.11 } \quad For each pair $ (p, q) $ above, by Chinese
Remainder Theorem (See [L]) and Dirichlet Theorem on primes in
arithmetic progressions [see[IK]], it is easy to verify that there
exist infinitely many $ D $ satisfying the conditions in every
Corollaries above.

\par \vskip 0.5 cm

\hspace{-0.6cm}{\bf 2 \ \ Proofs of Theorems}

\par \vskip 0.3 cm

Let $ M_{\Q} $ be the set of all places of $ \Q, $ including the
infinite $ \infty . $ For each place $ p, $ denote by $ \Q_{p} $ the
completion of $ \Q $ at $ p, $ and if $ p $ is finite, denote by $
v_{p} $ the corresponding normalized additive valuation, so $
v_{p}(p) = 1. $ Let $ S = \{ \infty , 2, p, q, D_{1}, \cdots , D_{n}
\}, $ and define a subgroup of $ \Q^{\star } / \Q^{\star ^{2}} $ as
follows: \ $ \Q(S, 2)= <-1, 2, p, q, D_{1},  \cdots , D_{n} > \cong
\left(\Z / 2 \Z \right)^{n + 4}. $ \ For any subset $ A \subset
\Q^{\star }, $ we write $ <A> $ for the subgroup of $ \Q^{\star } /
\Q^{\star ^{2}} $ generated by all the elements in $
A. $ For each $ d \in \Q(S, 2), $ define the curves \\
$ C_{d} : \ d w^{2} = d^{2} - 2 \varepsilon ( p + q) D d z^{2} + 4
D^{2} z^{4}, $ \ and \\
$ C_{d}^{\prime } : \ d w^{2} = d^{2} +  \varepsilon ( p + q) D d
z^{2}
+ p q D^{2} z^{4}. $ \\
Then by the theory of descent via two-isogeny ( see [S, chapter X]),
the $ \varphi -$Selmer group for the elliptic curve $ E / \Q
$ in (1.1) is \\
$ S^{(\varphi )} (E / \Q ) \cong \{ d \in \Q(S, 2) : \ C_{d} (
\Q_{v}) \neq \emptyset  \text{ for all } v \in S \}, $ \ and the $
\widehat{\varphi } -$Selmer group for the elliptic curve $
E^{\prime} / \Q $ in (1.2) is \\
$ S^{(\widehat{\varphi })} (E^{\prime } / \Q ) \cong \{ d \in
\Q(S, 2) : \ C_{d}^{\prime } ( \Q_{v}) \neq \emptyset  \text{ for
all } v \in S \}. $
\par  \vskip 0.2 cm

{\bf Proposition 2.1 } \ We assume $ \varepsilon = 1 $ and the
elliptic curve $ E = E_{+} $ be as in (1.1).
\par  \vskip 0.15 cm
(A) \ For $ d \in \Q (S, 2), $ if one of the following conditions holds: \\
(1) \ $ d < 0; $ \quad  (2) \ $ p \mid d; $ \quad  (3) \ $ q \mid
d. $ \\
Then $ d \notin S^{(\varphi )} (E / \Q ). $ \ Moreover, if $ d >
0, $ then $ C_{d}( \R ) \neq \emptyset . $
\par  \vskip 0.15 cm
(B) \ For $ d =2 \in \Q (S, 2), $ we have \\
(1) \ $ C_{2} ( \Q _{2}) \neq \emptyset \ \Longleftrightarrow \ D
( D - 2p - 2 ) \equiv 1 (\text{mod} 16 ). $ \\
(2) \ For each prime number $ l \mid p q D, \ C_{2} ( \Q _{l})
\neq \emptyset \ \Longleftrightarrow \ \left(\frac{2}{l}\right) =
1, $ i.e., $ l \equiv 1, 7 (\text{mod} 8 ). $
\par  \vskip 0.15 cm
(C) \ For each $ D_{i} \ ( 1\leq i \leq n ), $ we have \\
(1) \ $ C_{D_{i}} ( \Q _{2}) \neq \emptyset \ \Longleftrightarrow
\ D_{i} \equiv 1 (\text{mod} 4 ). $ \\
(2) \ For each prime number $ l \mid p q \widehat{D_{i}},  \
C_{D_{i}} ( \Q _{l}) \neq \emptyset \ \Longleftrightarrow \
\left(\frac{D_{i}}{l}\right) = 1. $  \\
(3) \ $ C_{D_{i}} ( \Q _{D_{i}}) \neq \emptyset \
\Longleftrightarrow \ \left(\frac{p \widehat{D_{i}}}{D_{i}}\right)
= \left(\frac{q \widehat{D_{i}}}{D_{i}}\right)= 1. $
\par  \vskip 0.2 cm

{\bf Proof. } \  Conclusion (A) follows directly by the valuation
property. \\
(B) \ In this case, $ C_{2}: w^{2} = 2 - 2 (p + q) D z^{2} + 2 D^{2}
z^{4}. $ \\
(1) \ If $ C_{2} ( \Q _{2}) \neq \emptyset , $ then there exists $ (
z_{0}, w_{0} ) \in \Q_{2}^{2} $ such that  $  w_{0}^{2} = 2 - 2 (p +
q ) D z_{0}^{2} + 2 D^{2} z_{0}^{4}, $ from which we have $
v_{2}(w_{0}) \geq 1 $ and $ v_{2}(z_{0}) = 0. $ Let $ w_{0} = 2
w_{1} $ with $ w_{1} \in \Q _{2} $ and $ v_{2}(w_{1}) \geq 0. $ Then
$ 2 w_{1}^{2} = 1 - (p + q ) D z_{0}^{2} + D^{2} z_{0}^{4}. $ Note
that $ v_{2}(D z_{0}) = 0, $ so by lemma in [Rob], we have $
z_{0}^{2} = 8 s + 1 $ and $ D^{2} z_{0}^{4} = 8 t + 1 $ for some $
s, t \in \Q _{2} $ with $ v_{2}(s) \geq 0 $ and $ v_{2}(t) \geq 0. $
So $ 2 w_{1}^{2} = 2 - 2 (p + 1) D (8 s + 1) + 8 t. $ Thus $ v_{2}(
2 w_{1}^{2} ) = v_{2}(2) = 1, $  so $ v_{2}(w_{1}) = 0, $ and we may
assume that $ (w_{1})^{2} = 8 r + 1 $ for some $ r \in \Q_{2} $ and
$ v_{2}(r) \geq 0. $ So  $ 2 ( 8 r + 1) = 1 - 2 (p + 1) D (8 s + 1)
+ D^{2} (8 s + 1 )^{2}, $ from which we get $ D^{2} - 2 (p + 1) D
\equiv 1 (\text{mod} 16 ), $ i.e., $ D
(D - 2 p -2) \equiv 1 (\text{mod} 16 ). $ \\
Conversely, if $ D (D - 2 p -2) \equiv 1 (\text{mod} 16 ), $  let $
f (z , w) = w^{2} - 2 + 2 (p + q) D z^{2} - 2 D^{2} z^{4}, $  then $
f'_{w} = 2 w, $  and $ v_{2}(f(1 , 2)) \geq 5 > 2 v_{2}(f'_{w}(1 ,
2)), $  so by Hensel's lemma (see [S],$ p. 322 $), we have $
C_{2}( \Q_{2}) \neq \emptyset . $ \\
(2) \ For $ l = p, $ if $ C_{2} ( \Q_{p} ) \neq \emptyset ,$ then $
\exists ( z_{0},w_{0} ) \in \Q_{p}^{2} $ s.t. $ w_{0}^{2} = 2 - 2 (p
+ q) D z_{0}^{2} + 2 D^{2} z_{0}^{4}. $ \ If $ v_{p}(z_{0}) > 0, $
then $ v_{p}(w_{0}^{2} - 2) \geq 1, $ so $ \left(\frac{2}{p}\right)
= 1. $ If $ v_{p}(z_{0}) = 0£¬$ then we can obtain that $ v_{p}(D
z_{0}^{2} - 1 ) = 0. $ Let $ w_{1} = w_{0}(D z_{0}^{2} -1)^{-1} \in
\Q_{p}^{*}, $ then $ v_{p}(w_{1}^{2} - 2) \geq 1, $ so $
\left(\frac{2}{p}\right) = 1. $ If $ v_{p}(z_{0}) < 0£¬$  then $
v_{p}(z_{0}^{-1}) > 0, $  so by the equality $
(\frac{w_{0}}{z_{0}^{2}})^{2} = 2 (z_{0}^{-1})^{4} - 4 (p + 1) D
z_{0}^{-2} + 2 D^{2} $ we get $ v_{p}(w_{1}) \geq 0 $ and $ v_{p} (
w_{1}^{2} - 2 D^{2} ) \geq 1 $ with $ w_{1} =
\frac{w_{0}}{z_{0}^{2}}. $ So $ \left(\frac{2 D^{2}}{p}\right) = 1,
$ and then $ \left(\frac{2}{p}\right) = 1. $ To sum up, we obtain $
\left(\frac{2}{p}\right) = 1. $ \\
Conversely, if $  \left(\frac{2}{p}\right) = 1, $  then $ 2 \equiv
a^{2}( \text{mod} p) $ for some $ a \in \Z $ with $ p \nmid a. $
Then for the above f(z, w), we have $ v_{p}(f(0, a)) \geq 1 > 0 = 2
v_{p}(f'_{w}(0, a)). $  So by Hensel's lemma , we have $ C_{2}(
\Q_{p})\neq \emptyset . $ \quad The case $ l = q, D_{i} \ ( 1 \leq i
\leq n ) $ are similar. This proves (B).
\par  \vskip 0.15 cm

(C) \ In this case, $ C_{D_{i}}:  w^{2} = D_{i} - 2 (p + q) D z^{2}
+ 4 D \widehat{D_{i}} z^{4}. $  \\
(1) \ If $ C_{D_{i}} ( \Q _{2}) \neq \emptyset , $ then $ \exists
(z_{0}, w_{0}) \in \Q_{2}^{2} $ s. t. $ w_{0}^{2} = D_{i} -
 2 (p + q) D z_{0}^{2} + 4 D \widehat{D_{i}} z_{0}^{4}. $ \\
If $ v_{2}(z_{0}) \geq 0, $ then $ v_{2}(w_{0}^{2}) = v_{2}(D_{i}),
\ v_{2}(w_{0}) = 0 $ and $ v_{2}(w_{0}^{2} - D_{i}) \geq 2. $ Since
$ v_{2}(w_{0}^{2} - D_{i}) \geq 3  $ ( see lemma in [Rob]), we have
$ v_{2}( D_{i} - 1 ) = v_{2}((w_{0}^{2} - 1) + ( D_{i} - w_{0}^{2}))
\geq 2. $  So $ D_{i} \equiv 1 (\text{mod} 4). $ \\
If $ v_{2}(z_{0}) < 0, $  then $ v_{2}(z_{0}^{-1}) \geq 1, $ so by $
(\frac{w_{0}}{2 z_{0}^{2}})^{2} = D_{i} ( 2 z_{0}^{2} )^{-2} - ( p +
1 ) D z_{0}^{-2} + D \widehat{D_{i}}, $ we get $
v_{2}((\frac{w_{0}}{2 z_{0}^{2}})^{2}) = v_{2}(D \widehat{D_{i}}) =
0, $ so  $ v_{2}((\frac{w_{0}}{2 z_{0}^{2}}) = 0 $ and $ v_{2}(
D_{i} \widehat{D_{i}}^{2} - (\frac{w_{0}}{2 z_{0}^{2}})^{2} ) \geq
2, $ hence $ v_{2}(D_{i} - (\frac{w_{0}}{2 z_{0}^{2}
\widehat{D_{i}}})^{2}) \geq 2. $ Since $ v_{2}((\frac{w_{0}}{2
z_{0}^{2} \widehat{D_{i}}})^{2} - 1 ) \geq 3, $  we get $
v_{2}(D_{i} - 1 ) = v_{2}((D_{i} - (\frac{w_{0}}{2 z_{0}^{2}
\widehat{D_{i}}})^{2}) + ( (\frac{w_{0}}{2 z_{0}^{2}
\widehat{D_{i}}})^{2} - 1 )) \geq 2. $  So $ D_{i} \equiv 1
(\text{mod} 4). $ \\
Conversely, if $ D_{i} \equiv 1 ( \text{mod} 4), $  then for the
polynomial $ f( z, w ) = w^{2} - D_{i} + 2 (p + q) D z^{2} - 4 D
\widehat{D_{i}} z^{4}, $  we have $ f'_{w} = 2 w. $ So if $ D_{i}
\equiv 1 (\text{mod} 8), $ then  $ v_{2}(f(0, 1)) \geq 3 > 2 = 2
v_{2}(f'_{w}(0, 1)); $  If $ D_{i} \equiv 5 (\text{mod} 8), $ then $
v_{2}(f(1, 1)) \geq 3 > 2 = 2 v_{2}(f'_{w}(1, 1)). $  So by Hensel's
lemma, we have $ C_{D_{i}} ( \Q _{2}) \neq \emptyset . $ \\
(2) \ For $ l = p, $ if $ C_{D_{i}} ( \Q _{p}) \neq \emptyset , $
then $ \exists (z_{0}, w_{0}) \in \Q_{p}^{2} $  s. t. $ w_{0}^{2} =
D_{i} - 2 (p + q) D z_{0}^{2} + 4 D \widehat{D_{i}} z_{0}^{4} = -4 p
D z_{0}^{2} + D_{i} (2 \widehat{D_{i}} z_{0}^{2} - 1)^{2}. $ \\
If $ v_{p}(z_{0}) > 0, $ then  $ v_{p}(w_{0}^{2}) = v_{p}(D_{i}) =
0, $ so $ v_{p}(w_{0}) = 0 $  and $  v_{p}(w_{0}^{2} - D_{i}) \geq
1, $  so $ \left(\frac{D_{i}}{p}\right) = 1. $ \\
If $ v_{p}(z_{0}) = 0, $  then $ v_{p}(w_{0}) \geq 0. $ If $ v_{p}(
2 \widehat{D_{i}} z_{0}^{2} - 1) \geq 1, $ then $ v_{p}(( 2
\widehat{D_{i}} z_{0}^{2} - 1))^{2} \geq 2 > 1 = v_{p}(-4 p D
z_{0}^{2}), $ so  $ v_{p}(w_{0}^{2}) = v_{p}(-4 p D z_{0}^{2}), $
i.e., $ 2 v_{p}(w_{0}^{2}) = 1, $ which is impossible !  So $ v_{p}(
2 \widehat{D_{i}} z_{0}^{2} - 1 ) = 0, $ and then we have $
v_{p}(D_{i} - (\frac{w_{0}}{2 \widehat{D_{i}} z_{0}^{2} - 1})^{2})
\geq 1, $ which shows that $ \left(\frac{D_{i}}{p}\right) = 1. $ \\
If $ v_{p}(z_{0}) < 0, $ then $ v_{p}(z_{0}^{-1}) \geq 1. $  By the
above equation of $ C_{D_{i}} $ we have $ (\frac{w_{0}}{2
\widehat{D_{i}} z_{0}^{2}})^{2} = D_{i} (2 \widehat{D_{i}}
z_{0}^{2})^{-2} - \frac{( p + 1 ) D}{(\widehat{D_{i}}
z_{0}^{2})^{2}} + D_{i}. $ So $ v_{p} ((\frac{w_{0}}{2
\widehat{D_{i}} z_{0}^{2}})^{2}) = v_{p}(D_{i}) = 0, $ hence $ v_{p}
(D_{i} - \frac{w_{0}}{2 \widehat{D_{i}} z_{0}^{2}} ) \geq 1, $ which
shows $ \left(\frac{D_{i}}{p}\right) = 1. $ \\
Conversely, if $ \left(\frac{D_{i}}{p}\right) = 1, $ then $ D_{i}
\equiv a^{2} (\text{mod} p) $ for some $ a \in \Z $  with $ v_{p}(a)
= 0. $ Then $ v_{p}(f(0, a)) \geq 1 > 0 = 2 v_{p}(f'_{w}(0, a)). $
By Hensel's lemma, we get $ C_{D_{i}} ( \Q _{p}) \neq \emptyset . $
The other cases follow the same line. We omit the details. The proof
of Proposition 2.1 is completed. \quad $ \Box $
\par  \vskip 0.2 cm

{\bf Corollary 2.2 } \ We assume $ \varepsilon = 1 $ and the
elliptic curve $ E = E_{+} $ be as in (1.1).
\par  \vskip 0.15 cm
(A) \ $ 2 \in S^{(\varphi )} (E / \Q ) \ \Longleftrightarrow \ p
\equiv 7 (\text{mod} 8 ) $ \ and \ $ D_{i} \equiv 1, 7 (\text{mod}
8 ) \ ( 1 \leq i \leq n ). $
\par  \vskip 0.15 cm
(B) \ For each $ D_{i} \ ( 1\leq i \leq n ), $ we have \ $ D_{i}
\in S^{(\varphi )} (E / \Q ) $ \ if and only if  \\ $ \ D_{i}
\equiv 1 (\text{mod} 4 ) $ \ and \ $
\left(\frac{D_{j}}{D_{i}}\right) = \left(\frac{p}{D_{i}}\right) =
\left(\frac{q}{D_{i}}\right) = 1 \quad ( 1 \leq j \leq n, \
\text{and } \ j \neq i). $
\par  \vskip 0.2 cm

{\bf Proof. } \ The results follow easily from Proposition 2.1.
\quad $ \Box$
\par  \vskip 0.3 cm

{\bf Proposition 2.3 } \ We assume $ \varepsilon = 1 $ and the
elliptic curve $ E^{\prime} = E^{\prime}_{+} $ be as in (1.2).
\par  \vskip 0.15 cm
(A) \ (1) \ For any $ d \in \Q (S, 2), \ C_{d}^{\prime } ( \R )
\neq \emptyset . $ If $ 2 \mid d, $ then $ d \notin S^{(\widehat{\varphi })}
(E^{\prime } / \Q ). $ \\
(2) \ $ \{ 1, pq, -pD, -qD \} \subset S^{(\widehat{\varphi })}
(E^{\prime } / \Q ). $
\par  \vskip 0.15 cm
(B) \ For each $ D_{i} \ ( 1\leq i \leq n ), $ we have \\
(1) \ $ C_{D_{i}}^{\prime} ( \Q _{2}) \neq \emptyset $ if and only
if one of the following conditions holds: \\
(a) \ $ D_{i} \equiv
1 (\text{mod} 8 ); $ \quad (b) \ $ ( 1 + p \widehat{D_{i}} ) ( 1 +
q \widehat{D_{i}} ) \equiv 0 (\text{mod} 16 ); $ \\
(c) \ $ D_{i} \equiv 3 (\text{mod} 8 ) $ and $ p \equiv 1
(\text{mod} 4 ); $ \quad (d) \ $ D_{i} \equiv 7 (\text{mod} 8 ) $
and $ p \equiv 3 (\text{mod} 4 ). $ \\
(2) \ $ C_{D_{i}}^{\prime} ( \Q _{p}) \neq \emptyset $ and $
C_{D_{i}}^{\prime} ( \Q _{q}) \neq \emptyset . $ \\
(3) \ For each $ j \neq i, \quad  C_{D_{i}}^{\prime} ( \Q
_{D_{j}}) \neq \emptyset \Longleftrightarrow \ \left( 1 -
\left(\frac{D_{i}}{D_{j}}\right) \right) \left( 1 - \left(\frac{p
q D_{i}}{D_{j}}\right) \right) = 0. $ \\
(4) \ $ C_{D_{i}}^{\prime} ( \Q _{D_{i}}) \neq \emptyset
\Longleftrightarrow \ \left( 1 - \left( \frac{ - p
\widehat{D_{i}}}{D_{i}}\right) \right) \left( 1 - \left(\frac{- q
\widehat{D_{i}}}{D_{i}}\right) \right) = 0. $
\par  \vskip 0.15 cm
(C) \ For $ d = -pq \in \Q (S, 2), $ we have \\
(1) \ $ C_{-pq}^{\prime} ( \Q _{2}) \neq \emptyset \
\Longleftrightarrow \ p \equiv 3 (\text{mod} 4 ) $ \ or \ $
D - p \equiv 0, 2 (\text{mod} 8 ). $ \\
(2) \ $ C_{-pq}^{\prime} ( \Q _{p}) \neq \emptyset $ and $
C_{-pq}^{\prime} ( \Q _{q}) \neq \emptyset . $ \\
(3) \ For each $ D_{i} \ ( 1\leq i \leq n ), $ we have \\
$ C_{-pq}^{\prime} ( \Q _{D_{i}}) \neq \emptyset
\Longleftrightarrow \ \left( 1 - \left(\frac{-1}{D_{i}}\right)
\right) \left( 1 - \left(\frac{ - p q }{D_{i}}\right) \right) = 0.
$
\par  \vskip 0.2 cm

{\bf Proof. } \ (A) \ (1) \ Obvious. \\
(2) \ For $ C'_{d}: d w^{2} = (d + p D z^{2})(d + q D z^{2}), $ if $
d = - p D, $ then $ (1, 0) \in C'_{d}(\Q); $ and  if $ d = - q D, $
then $ (1, 0) \in C'_{d}(\Q), $ so the group $ < - p D, - q D
> =  \{ 1, pq, -pD, -qD \} \subseteq
E(\Q)/ \widehat{\varphi }(E^{\prime }(\Q )) \subset
S^{(\widehat{\varphi })} (E^{\prime } / \Q ). $  \\
(B) \ $ C_{D_{i}}^{\prime } : \ w^{2} = D_{i} + (p + q) D z^{2} + p
q D \widehat{D_{i}} z^{4} = D_{i}(1 + p \widehat{D_{i}} z^{2})(1 + q
\widehat{D_{i}} z^{2}). $   \\
(1) \ If $  C_{D_{i}}^{\prime }( \Q_{2}) \neq \emptyset , $ then $
\exists (z_{0}, w_{0}) \in \Q_{D_{i}}^{2} $  s. t. $ w_{0}^{2} =
D_{i} + (p + q) D z_{0}^{2} + p q D \widehat{D_{i}} z_{0}^{4}. $ \\
If $ v_{2}(z_{0}) > 0, $  then $ v_{2}(w_{0}) = 0. $ By Lemma 4 in
[C, p.11], we get $ v_{2}(D_{i} - 1) \geq 3 $ and then $ D_{i}
\equiv 1 (\text{mod} 8). $  \\
If $ v_{2}(z_{0}) \leq 0, $ then by the equality $
(\frac{w_{0}}{z_{0}^{2}})^{2} = D_{i}(z_{0}^{-2} + p
\widehat{D_{i}})(z_{0}^{-2} + q \widehat{D_{i}}), $ we have $
v_{2}(\frac{w_{0}}{z_{0}^{2}}) \geq 0. $ Denote $ W =
\frac{w_{0}}{z_{0}^{2}}, \ Z = z_{0}^{-1}, $  then  $ v_{2}(Z) \geq
0, \ v_{2}(W) \geq 0 $ and $ W^{2} = D_{i}(Z^{2} + p
\widehat{D_{i}})(Z^{2} + q \widehat{D_{i}}) = D_{i}(V^{2} -
\widehat{D_i}^2) $ with $ V = Z^{2} + (p + 1) \widehat{D_{i}} \in
\Q_{2}. $ Obviously, $ v_{2}(V) \geq 0. $ \\
i) \ If $ v_{2}(V) = 0, $ then  $ v_{2}(V^{2} - 1) \geq 3 $ and $
v_{2} (Z) = 0. $ So $ v_{2}(V^{2} - \widehat{D_{i}}^2) = v_{2}( (
V^{2} - 1 ) - ( \widehat{D_{i}}^2 - 1 ) ) \geq 3. $ Thus $
v_{2}(W^{2}) = v_{2}(V^{2} - \widehat{D_i}^{2}) \geq 3, $ so $
v_{2}(W^{2}) \geq 2, $ and then $ v_{2}((1 + p \widehat{D_{i}})(1 +
q \widehat{D_{i}})) \geq 4, $ i.e.,  $ ( 1 + p \widehat{D_{i}} ) ( 1
+ q \widehat{D_{i}} ) \equiv 0 (\text{mod} 16 ). $ \\
ii) \ If $ v_{2}(V) \geq 1, $ then  $ v_{2}(W) = 0. $ \\
If $ v_{2}(V) = 1, $ then $ v_{2}(V^{2} - 4) \geq 5. $ So $
v_{2}(W^{2} - 3 D_{i}) = v_{2}((V^{2} - 4) - ( \widehat{D_{i}}^{2} -
1)) \geq 3. $ Thus $ v_{2}(1 - 3 D_{i} ) \geq 3, $ hence $ D_{i}
\equiv 3 (\text{mod} 8). $ Since $ V = Z^{2} + (p + 1)
\widehat{D_i}, $ we get $ v_{2}(p + 1) = v_{2}(V - Z^{2}) = 1, $ so
$ p \equiv 1 (\text{mod} 4). $ \\
If $ v_{2}(V) \geq 2, $ then $ v_{2}(V^2) \geq 4, $ so $ W^{2} +
D_{i} = D_{i}(V^{2} + 1 - \widehat{D_{i}}^2 ), $ hence $ v_{2}(W^{2}
+ D_{i} ) \geq 3, $ and then $ v_{2}(1 + D_{i} ) = v_{2}(W^{2} +
D_{i} + 1 - W^{2} ) \geq 3, $ so $ D_{i} \equiv 7 (\text{mod} 8 ). $
Since $ V = Z^{2} + (p + 1) \widehat{D_{i}}, $ we get $ v_{2}(p + 1)
= v_{2}(V - Z^{2} ) \geq 2, $ which shows $ p \equiv 3 (\text{mod}
4). $ \\
To sum up, we have proved that if \ $ C_{D_{i}}^{\prime } ( \Q _{2})
\neq \emptyset , $ then one of the following conditions holds: \\
(a) \ $ D_{i} \equiv 1 (\text{mod} 8 ); $ \quad (b) \ $ ( 1 + p
\widehat{D_{i}} ) ( 1 + q \widehat{D_{i}} ) \equiv 0 (\text{mod} 16 ); $ \\
(c) \ $ D_{i} \equiv 3 (\text{mod} 8 ) $ and $ p \equiv 1
(\text{mod} 4 ); $ \quad (d) \ $ D_{i} \equiv 7 (\text{mod} 8 ) $
and $ p \equiv 3 (\text{mod} 4 ). $ \\
Conversely, we assume one of the above conditions holds. Let $ f(z,
\ w) =  D_{i} + (p + q) D z^{2} + p q D \widehat{D_{i}} z^{4} -
w^{2}, $ then $ f'_{w}(z, w) = -2 w. $ \\
(a) \ If $ D_{i} \equiv 1 (\text{mod} 8 ), $ then $ v_{2}(f(0, 1))
\geq 3 > 2 = 2 v_{2}(f'_{w}(0, 1)). $ So by Hensel's
lemma, we get $ C_{D_{i}}^{\prime } ( \Q _{2}) \neq \emptyset . $ \\
(b) \ If \ $ ( 1 + p \widehat{D_{i}} ) ( 1 + q \widehat{D_{i}} )
\equiv 0 (\text{mod} 16 ), $ then $ v_{2}(1 + p \widehat{D_{i}})
\geq 3 $ or $ v_{2}(1 + q \widehat{D_{i}}) \geq 3. $ Let $ g(z) = 1
+ l \widehat{D_{i}} z^{2}, $  where $ l = p $ or $ q $ according to
the above first or second inequalities. Then $ g'(z) = 2 l
\widehat{D_{i}} z. $ So $ v_{2}(g(1)) \geq 3 > 2 = 2 v_{2}(g'(1)), $
and by Hensel's lemma, $ g(z_{0}) = 0 $ for some $ z_{0} \in \Q_2. $
Therefore $ (z_{0}, 0) \in C_{D_{i}}^{\prime } ( \Q
_{2}). $ \\
(c) \ If $ D_{i} \equiv 3 (\text{mod} 8 ) $ and $ p \equiv 1
(\text{mod} 4 ), $ let $ g(w) = w^2 - D_{i}(2^2 + p
\widehat{D_{i}})(2^2 + q \widehat{D_{i}}), $ then $ g'(w) = 2 w, $
and $ v_{2}(g(1)) > 2 = 2 v_{2}(g'(1)). $ So by Hensel's lemma, $
g(w_{0}) = 0 $ for some $ w_{0} \in \Q_2. $ Thus $ (\frac{1}{2}, \
\frac{w_{0}}{4}) \in C_{D_{i}}^{\prime } ( \Q_{2}). $ \\
(d) \ similar to (c) above. This proves (1).  \\
(2) \ Let $ f(z) = D_{i} + (p + q) D z^{2} + p q D \widehat{D_{i}}
z^{4}, $ then  $ f(z) \equiv D_{i} (1 + 2 \widehat{D_{i}} z^{2})
(\text{mod} p ), $ also $ f(z) \equiv D_{i} (1 - 2 \widehat{D_{i}}
z^{2}) (\text{mod} q). $ \ Take $ \overline{g} = 1, \ \overline{h} =
\overline{f} (\text{mod} p ) $ \ ( if $ D_{i} $ is a square $
\text{mod} p, $ then take $ \overline{g} = D_{i} (\text{mod} p ) ),
$ then $ \text{deg} \overline{g} = 0 $ and it is easy to verify that
$ \overline{h} $ is square free. So by lemma 14 in [MSS], $ w^{2} =
f(z) $ has solutions in $ \Q_{p}, $ i.e., $ C_{D_{i}}^{\prime } ( \Q
_{p}) \neq \emptyset . $ Similarly, $ C_{D_{i}}^{\prime } ( \Q _{q})
\neq \emptyset . $ This proves (2).  \\
(3) \ If $ C_{D_{i}}^{\prime }(\Q_{D_{j}}) \neq \emptyset \ ( j \neq
i), $ then $ \exists (z_{0}, w_{0}) \in \Q_{D_{j}}^{2} $  s. t. $
w_{0}^{2} = D_{i}(1 + p \widehat{D_{i}} z_{0}^{2} ) ( 1 + q
\widehat{D_{i}} z_{0}^{2}). $  \\
i) \ If $ v_{D_{j}}(z_{0}) \geq 0, $ then $ v_{D_{j}}(w_{0}) = 0 $
and $ D_{i} \equiv w_{0}^{2} (\text{mod} D_{j} ), $ so $
\left(\frac{D_{i}}{D_{j}}\right) = 1. $ \\
ii) \ If $ v_{D_{j}}(z_{0}) < 0, $ then $ v_{D_{j}}(z_{0}^{-1}) \geq
1 $ and $ v_{D_{j}}(w_0) = 1 + 2 v_{D_{j}}(z_{0}) < 0. $ Now $
(\frac{w_{0}}{z_{0}^{2}})^{2} = D_{i}z_{0}^{-4} + (p + q) D
z_{0}^{-2} + p q D_{i} \widehat{D_{i}}^{2}. $ By the above
discussion, we get $ p q D_{i} \equiv (\frac{w_{0}}{z_{0}^{2}
\widehat{D_{i}}})^{2} ( \text{mod} D_{j} ) $ with $ v_{D_{j}}
(\frac{w_{0}}{z_{0}^{2}\widehat{D_{i}}}), $ so $ \left(\frac{p q
D_{i}}{D_{j}} \right) = 1. $ To sum up, we get $ \left( 1 -
\left(\frac{D_{i}}{D_{j}} \right) \right) \left( 1 - \left(\frac{p q
D_{i}}{D_{j}} \right) \right) = 0. $ \\
Conversely, if $ \left( 1 - \left(\frac{D_{i}}{D_{j}} \right)
\right) \left( 1 - \left(\frac{p q D_{i}}{D_{j}} \right) \right) =
0, $ then $ D_{i} \equiv a^{2} ( \text{mod} D_{j} ) $ or $ p q D_{i}
\equiv b^{2} ( \text{mod} D_{j}) $ for some $ a, b \in \Z $ and $
D_{j} \nmid a b. $ \\
For the former, let $ f(z , w) = w^{2} - D_{i} - (p + q) D z^{2} + p
q D_{i} \widehat{D_{i}}^{2} z^{4}, $ then $ f'_{w}(z, w) = 2 w $ and
$ v_{D_{j}}(f(0, a)) > 0 = 2 v_{D_{j}}(f'_{w}(0, a)). $ So by
Hensel's lemma, $ C_{D_{i}}^{\prime } ( \Q _{D_{j}}) \neq
\emptyset . $  \\
For the later, let $ g(z, w) = w^{2} - D_{i} \widehat{D_{i}}^{2} -
(p + q) \widehat{D_{i}}^{2} D z^{2} - p q D_{i} \widehat{D_{i}}^{4}
z^{4}, $ then $ g'_{w}(z, w) = 2 w $ and $ v_{D_{j}}( g(
\widehat{D_{i}}^{-1}, b)) > 0 = 2 v_{D_{j}}(g'_{w}(
\widehat{D_{i}}^{-1}, b)). $ So by Hensel's lemma, $ \exists (z_{0},
w_{0}) \in  \Q _{D_{j}}^{2} $ s. t. $ g(z_{0}, w_{0}) = 0, $ and so
$ ( z_{0}, \frac{w_{0}}{\widehat{D_{i}}}) \in C_{D_{i}}^{\prime }
(\Q_{D_{j}}). $ This proves (3). \\
(4) \ If $ C_{D_{i}}^{\prime }(\Q_{D_{i}}) \neq \emptyset , $ then $
\exists (z_{0}, w_{0}) \in \Q_{D_{i}}^{2} $ s. t. $  w_{0}^{2} =
D_{i}(1 + p \widehat{D_{i}} z_{0}^{2})(1 + q \widehat{D_{i}}
z_{0}^{2}). $ It is easy to see that $ v_{D_{i}}( z_{0}) = 0 $ and $
v_{D_{i}}( w_{0}) > 0. $ Set $ w_{0} = D_{i} w_{1} $  with $ w_{1}
\in \Q_{D_{i}} $ and $ v_{D_{i}} ( w_{1}) \geq 0. $ Then $ D_{i}
w_{1}^{2} = ( 1 + p \widehat{D_{i}} z_{0}^{2})(1 + q \widehat{D_{i}}
z_{0}^{2}), $ so $ ( 1 + p \widehat{D_{i}} z_{0}^{2})(1 + q
\widehat{D_{i}} z_{0}^{2}) \equiv 0 (\text{mod} D_{i}). $ Thus $
\left( 1 - \left( \frac{ - p \widehat{D_{i}}}{D_{i}} \right) \right)
\left( 1 - \left(\frac{- q \widehat{D_{i}}}{D_{i}} \right) \right) = 0. $ \\
Conversely, if $ \left( 1 - \left( \frac{ - p
\widehat{D_{i}}}{D_{i}} \right) \right) \left( 1 - \left(\frac{- q
\widehat{D_{i}}}{D_{i}} \right) \right) = 0, $ then $ \exists a \in
\Z $  such that $ -l \widehat{D_{i}} \equiv a^{2} (\text{mod} D_{i})
$ for $ l = p $ or $ q. $ Obviously, $ D_{i} \nmid a. $ \ Let $ g(z)
= 1 + l \widehat{D_{i}} z^{2}, $ then $ g'(z) = 2 l \widehat{D_{i}}
z. $ So $ v_{D_{i}}( g(a^{-1})) > 0 = 2 v_{D_{i}}( g'(a^{-1})), $ by
Hensel's lemma, $ g(z_{0}) = 0 $ for some $ z_{0} \in \Q_{D_{i}}. $
Thus $ ( z_{0} , 0) \in C_{D_{i}}^{\prime }(\Q_{D_{i}}), $ so $
C_{D_{i}}^{\prime }(\Q_{D_{i}}) \neq \emptyset . $ This proves (4),
and (B) is proved.  \\
(C) \ (1) \ If $ C_{-pq}^{\prime }(\Q_{2}) \neq \emptyset , $ then $
\exists (z_{0}, w_{0}) \in \Q_{2}^{2} $ s. t. $ w_{0}^{2} =
- p q + (p + q) D z_{0}^{2} - D^{2} z_{0}^{4}. $ \\
i) \ If $ v_{2}(z_{0}) < 0, $ then $ v_{2}(z_{0}^{-1}) \geq 1, $ so
by the equality $ r_{0}^{2} = - p q z_{0}^{-4} + (p + q)D z_{0}^{-2}
- D^{2} $ with $ r_{0} = \frac{w_{0}}{z_{0}^{2}}, $ we get $
v_{2}(r_{0}) = 0, \ v_{2}(r_{0}^{2} + D^{2} ) \geq 4. $ Thus $
v_{2}( r_{0}^{2} + D^{2}) = v_{2}((r_{0}^{2} - D^{2} ) + 2 D^{2} ) =
1, $ a contradiction! So we must have $ v_{2}(z_{0}) \geq 0. $ \\
ii) \ If $ v_{2}(z_{0}) > 0, $ then $ v_{2}(w_{0}) = 0, \ v_{2}(
w_{0}^{2} + p q ) \geq 4. $ So we get $ v_{2}(1 + p q) = v_{2}(1 -
w_{0}^{2} +  w_{0}^{2} + p q ) \geq 3, $  which shows $ p
\equiv 3 ( \text{mod} 4). $ \\
iii) \ If $ v_{2}(z_{0}) = 0, $ then $ v_{2}(z_{0}^{2} - 1) \geq 3,
\ v_{2}(z_{0}^{4} - 1) \geq 4. $ Since $ w_{0}^{2} = - (D - p)(D -
q) + (p + q) D (z_{0}^{2} -1 ) - D^{2} (z_{0}^{4} - 1 ), $ we have $
v_{2}( w_{0}^{2} + (D - p)(D - q) ) \geq 4. $ But $ (D - p)(D - q)
\equiv 0 ( \text{mod} 8 ), $ so $ v_{2}(w_{0}) \geq 2. $ Thus $
v_{2}((D - p)(D - q)) \geq 4, $  which shows $ D - p \equiv 0, \ 2
(\text{mod} 8). $ \\
Conversely, if $ p \equiv 3 ( \text{mod} 4), $ we get $ -p q \equiv
1 (\text{mod} 8). $ So $ -p q = a^{2} $ for some $ a \in \Q_{2}. $
Hence $ (0 , a) \in C_{-pq}^{\prime }(\Q_{2}). $ \\
If $ D - p \equiv 0, \ 2 (\text{mod} 8 ), $ then $ (D - p - 1)^{2}
\equiv 1 (\text{mod} 16). $ Let $ f(z, w)= - [(p + 1) z^{2} - D]^{2}
+ z^{4} - w^{2}, $ then $ f'_{z}(z, w)= -4 (p + 1) [(p + 1) z^{2} -
D] z + 4 z^{3}. $ \\
If $ v_{2}((D - p - 1)^{2} - 1) > 4, $ then  $ v_{2}( f(1, 0)) > 4 =
2 v_{2}( f'_{z}(1, 0)); $ \\
If $ v_{2}((D - p - 1)^{2} - 1) = 4, $ then  $ v_{2}( f(3, 0)) \geq
5 > 4 = 2 v_{2}( f'_{z}(3, 0)). $ So by Hensel's lemma, $ \exists
(z_{0}, w_{0}) \in \Q_{2}^{2} $ s. t. $ f(z_{0}, w_{0}) = 0. $ It is
easy to prove $ z_{0} \neq 0. $ Hence $ (z_{0}^{-1}, w_{0} z_{0}^{-1})
\in C_{-pq}^{\prime }(\Q_{2}). $ \\
(2) is similar to that of (B)(2). \\
(3) \ If $ C_{-pq}^{\prime }(\Q_{D_{i}}) \neq \emptyset , $ then $
\exists (z_{0}, w_{0}) \in \Q_{D_{i}}^{2} $ s. t. $ w_{0}^{2} = - p
q + (p + q)D z_{0}^{2} - D^{2} z_{0}^{4}. $ \\
i) \ If $ v_{D_{i}}(z_{0}) \geq 0, $ then we get
$ \left(\frac{ -p q }{D_{i}} \right) = 1. $ \\
ii) \ If $ v_{D_{i}}(z_{0}) < 0, $ then $ v_{D_{i}}(z_{0}^{-1}) \geq
1. $ Since $ (\frac{w_{0}}{D z_{0}^{2}})^{2} = -D^{-2} p q
z_{0}^{-4} + (p + q) D^{-1} z_{0}^{-2} - 1, $ we have $ \left(\frac{
-1 }{D_{i}} \right) = 1. $  To sum up, we get $ \left( 1 -
\left(\frac{-1}{D_{i}}\right) \right) \left( 1 - \left(\frac{ - p q
}{D_{i}}\right) \right) = 0. $ \\
Conversely, if $ \left( 1 - \left(\frac{-1}{D_{i}}\right) \right)
\left( 1 - \left(\frac{ - p q }{D_{i}}\right) \right) = 0, $ then $
\exists a, b \in \Z $ s. t. $ -p q \equiv a^{2} (\text{mod} D_{i}) $
or $ -1 \equiv b^{2} (\text{mod} D_{i}) $ with $ D_{i} \nmid a b. $
For the former, let $ f(z, w) = w^{2} + p q - (p + q)  D z^{2} + D^2
z^{4}, $ then $ f'_{w}(z, w) = 2 w $  and  $ v_{D_{i}}(f(0, a)) > 0
= 2 v_{D_{i}}( f'_{w}(0, a)). $ So by  Hensel's lemma, we get $
C_{-pq}^{\prime } ( \Q_{D_{i}}) \neq \emptyset . $ \\
For the later, let $ f(z, w) = w^{2} + D^{2} p q - (p + q) D^{3}
z^{2} + D^{4} z^{4}, $ then $ f'_{w}(z, w) = 2 w $ and $ v_{D_{i}}(
f(D^{-1}, b)) > 0 = 2 v_{D_{i}}(f'_{w}(D^{-1}, b)). $ So by Hensel's
lemma, $ \exists (z_{0}, w_{0}) \in \Q_{D_{i}}^{2} $ s. t. $
f(z_{0}, w_{0}) = 0. $  Thus $ \exists (z_{0}, w_{0} D^{-1}) \in
C_{-pq}^{\prime } ( \Q_{D_{i}}). $ This proves (C), and the proof of
proposition 2.3 is completed. \quad $ \Box $
\par  \vskip 0.2 cm

{\bf Corollary 2.4 } \ We assume $ \varepsilon = 1 $ and the
elliptic curve $ E^{\prime} = E^{\prime}_{+} $ be as in (1.2).
\par  \vskip 0.15 cm
(A) Denote \ $ \alpha ( -pq )^{\prime} = \sum_{i = 1}^{n} \left( 1
- \left(\frac{-1}{D_{i}}\right) \right) \left( 1 - \left(\frac{ -
p q }{D_{i}}\right) \right). $ \ Then \\ $ -pq \in
S^{(\widehat{\varphi })} (E^{\prime } / \Q ) \ \Longleftrightarrow
\ \left \{
   \begin{array}{l}
  p \equiv 3 (\text{mod} 4 ), \\
  \alpha ( -pq )^{\prime} = 0;
  \end{array}
  \right.   \ \text{or} \ \left \{ \begin{array}{l}
  D - p \equiv 0, 2 (\text{mod} 8 ), \\
  \alpha ( -pq )^{\prime} = 0.
  \end{array}
  \right. $
\par  \vskip 0.15 cm
(B) \ Denote \ $ \beta ( -D )^{\prime} = \sum_{i = 1}^{n} \left( 1
- \left(\frac{p}{D_{i}}\right) \right) \left( 1 - \left(\frac{q
}{D_{i}}\right) \right). $ \ Then \ $ -D \in S^{(\widehat{\varphi
})} (E^{\prime } / \Q ) $ \ if
and only if one of the following conditions holds: \\
(1) \ $ \left \{
   \begin{array}{l}
  D \equiv 7 (\text{mod} 8 ), \\
  \beta ( -D )^{\prime} = 0;
  \end{array}
  \right.   \qquad \quad (2) \ \left \{ \begin{array}{l}
  p \equiv 1, 7 (\text{mod} 8 ), \\
  \beta ( -D )^{\prime} = 0;
  \end{array}
  \right. $ \\
(3) \ $ \left \{
   \begin{array}{l}
  D \equiv 1 (\text{mod} 8 ), \\
  p \equiv 3 (\text{mod} 8 ), \\
  \beta ( -D )^{\prime} = 0;
  \end{array}
  \right.   \qquad  \quad (4) \ \left \{ \begin{array}{l}
  D \equiv 5 (\text{mod} 8 ), \\
  p \equiv 5 (\text{mod} 8 ), \\
  \beta ( -D )^{\prime} = 0.
  \end{array}
  \right. $
\par  \vskip 0.15 cm
(C) \ For each $ D_{i} \ ( 1\leq i \leq n ), $ denote $$ \Pi
 _{i}^{+}(D)^{\prime } = \left( 1 - \left(\frac{ - p
\widehat{D_{i}}}{D_{i}}\right) \right) \left( 1 - \left(\frac{- q
\widehat{D_{i}} }{D_{i}}\right) \right) + \sum_{j = 1, \ j \neq
i}^{n} \left( 1 - \left(\frac{D_{i}}{D_{j}}\right) \right) \left(
1 - \left(\frac{p q D_{i} }{D_{j}}\right) \right). $$ \ Then \ $
D_{i} \in S^{(\widehat{\varphi })} (E^{\prime } / \Q ) $ \ if
and only if one of the following conditions holds: \\
(1) \ $ \left \{
   \begin{array}{l}
  D_{i} \equiv 1 (\text{mod} 8 ), \\
  \Pi _{i}^{+}(D)^{\prime } = 0;
  \end{array}
  \right.   \qquad \quad (2) \ \left \{ \begin{array}{l}
  (1 + p \widehat{D_{i}} )(1 + q \widehat{D_{i}} )  \equiv 0 (\text{mod} 16 ), \\
  \Pi _{i}^{+}(D)^{\prime } = 0;
  \end{array}
  \right. $ \\
(3) \ $ \left \{
   \begin{array}{l}
  D_{i} \equiv 3 (\text{mod} 8 ), \\
  p \equiv 1 (\text{mod} 4 ), \\
  \Pi _{i}^{+}(D)^{\prime } = 0;
  \end{array}
  \right.   \qquad  \quad (4) \ \left \{ \begin{array}{l}
  D_{i} \equiv 7 (\text{mod} 8 ), \\
  p \equiv 3 (\text{mod} 4 ), \\
  \Pi _{i}^{+}(D)^{\prime } = 0.
  \end{array}
  \right. $
\par  \vskip 0.2 cm

{\bf Proof. } \ The results follow easily from Proposition 2.3.
\quad $ \Box$
\par  \vskip 0.2 cm

{\bf Proposition 2.5 } \ We assume $ \varepsilon = - 1 $ and the
elliptic curve $ E = E_{-} $ be as in (1.1).
\par  \vskip 0.15 cm
(A) \ For $ d \in \Q (S, 2), $ if one of the following conditions holds: \\
(1) \ $ p \mid d; $ \quad  (2) \ $ q \mid d; $ \quad (3) \ $ d = -1. $\\
Then $ d \notin S^{(\varphi )} (E / \Q ). $
\par  \vskip 0.15 cm
(B$_{1}$) \ For $ d = 2 \in \Q (S, 2), $ we have \\
(1) \ $ C_{2} ( \Q _{2}) \neq \emptyset \ \Longleftrightarrow \ D
( D + 2p + 2 ) \equiv 1 (\text{mod} 16 ). $ \\
(2) \ For each prime number $ l \mid p q D, \ C_{2} ( \Q _{l})
\neq \emptyset \ \Longleftrightarrow \ \left(\frac{2}{l}\right) =
1, $ i.e., $ l \equiv 1, 7 (\text{mod} 8 ). $
\par  \vskip 0.15 cm
(B$_{2}$) \ For $ d = - 2 \in \Q (S, 2), $ we have \\
(1) \ $ C_{- 2} ( \Q _{2}) \neq \emptyset \ \Longleftrightarrow \
D( - D + 2p + 2 ) \equiv 3 (\text{mod} 16 ). $ \\
(2) \ For each prime number $ l \mid p q D, \ C_{ - 2} ( \Q _{l})
\neq \emptyset \ \Longleftrightarrow \ \left(\frac{- 2}{l}\right)
= 1, $ i.e., $ l \equiv 1, 3 (\text{mod} 8 ). $
\par  \vskip 0.15 cm
(C$_{1}$) \ For each $ D_{i} \ ( 1\leq i \leq n ), $ we have \\
(1) \ $ C_{D_{i}} ( \Q _{2}) \neq \emptyset \ \Longleftrightarrow
\ D_{i} \equiv 1 (\text{mod} 4 ). $ \\
(2) \ For each prime number $ l \mid p q \widehat{D_{i}},  \
C_{D_{i}} ( \Q _{l}) \neq \emptyset \ \Longleftrightarrow \
\left(\frac{D_{i}}{l}\right) = 1. $  \\
(3) \ $ C_{D_{i}} ( \Q _{D_{i}}) \neq \emptyset \
\Longleftrightarrow \ \left(\frac{- p
\widehat{D_{i}}}{D_{i}}\right) = \left(\frac{ - q
\widehat{D_{i}}}{D_{i}}\right)= 1. $
\par  \vskip 0.15 cm
(C$_{2}$) \ For each $ D_{i} \ ( 1\leq i \leq n ), $ we have \\
(1) \ $ C_{ - D_{i}} ( \Q _{2}) \neq \emptyset \
\Longleftrightarrow
\ D_{i} \equiv 3 (\text{mod} 4 ). $ \\
(2) \ For each prime number $ l \mid p q \widehat{D_{i}},  \ C_{-
D_{i}} ( \Q _{l}) \neq \emptyset \ \Longleftrightarrow \
\left(\frac{ - D_{i}}{l}\right) = 1. $  \\
(3) \ $ C_{- D_{i}} ( \Q _{D_{i}}) \neq \emptyset \
\Longleftrightarrow \ \left(\frac{ p
\widehat{D_{i}}}{D_{i}}\right) = \left(\frac{ q
\widehat{D_{i}}}{D_{i}}\right)= 1. $

\par  \vskip 0.2 cm

{\bf Proof. } \ The proof is similar to Proposition 2.1, we omit the
details.

\par  \vskip 0.2 cm

{\bf Corollary 2.6 } \ We assume $ \varepsilon = - 1 $ and the
elliptic curve $ E = E_{-} $ be as in (1.1).
\par  \vskip 0.15 cm
(A$_{1}$) \ $ 2 \in S^{(\varphi )} (E / \Q ) \ \Longleftrightarrow
\ p \equiv 7 (\text{mod} 8 ) $ \ and \ $ D_{i} \equiv 1, 7
(\text{mod} 8 ) \ ( 1 \leq i \leq n ). $
\par  \vskip 0.15 cm
(A$_{2}$) \ $ - 2 \in S^{(\varphi )} (E / \Q ) \
\Longleftrightarrow \ p \equiv 1 (\text{mod} 8 ) $ \ and \ $ D_{i}
\equiv 1, 3 (\text{mod} 8 ) \ ( 1 \leq i \leq n ). $
\par  \vskip 0.15 cm
(B$_{1}$) \ For each $ D_{i} \ ( 1\leq i \leq n ), $ we have \ $
D_{i} \in S^{(\varphi )} (E / \Q ) $ \ if and only if  \\ $ \
D_{i} \equiv 1 (\text{mod} 4 ) $ \ and \ $
\left(\frac{D_{j}}{D_{i}}\right) = \left(\frac{p}{D_{i}}\right) =
\left(\frac{q}{D_{i}}\right) = 1 \quad ( 1 \leq j \leq n, \
\text{and } \ j \neq i). $
\par  \vskip 0.15 cm
(B$_{2}$) \ For each $ D_{i} \ ( 1\leq i \leq n ), $ we have \ $ -
D_{i} \in S^{(\varphi )} (E / \Q ) $ \ if and only if  \\ $ \
D_{i} \equiv 3 (\text{mod} 4 ) $ \ and \ $
\left(\frac{D_{j}}{D_{i}}\right) = \left(\frac{p}{D_{i}}\right) =
\left(\frac{q}{D_{i}}\right) = 1 \quad ( 1 \leq j \leq n, \
\text{and } \ j \neq i). $
\par  \vskip 0.2 cm

{\bf Proof. } \ The results follow easily from Proposition 2.5.
\quad $ \Box$
\par  \vskip 0.2 cm

{\bf Proposition 2.7 } \ We assume $ \varepsilon = - 1 $ and the
elliptic curve $ E^{\prime} = E^{\prime}_{-} $ be as in (1.2).
\par  \vskip 0.15 cm
(A) \ (1) \ For any $ d \in \Q (S, 2) $ and $ d > 0, \ C_{d}^{\prime
} ( \R ) \neq \emptyset . $ If $ 2 \mid d $ or $ d < 0, $ then $ d
\notin S^{(\widehat{\varphi })}
(E^{\prime } / \Q ). $ \\
(2) \ $ \{ 1, pq, pD, qD \} \subset S^{(\widehat{\varphi })}
(E^{\prime } / \Q ). $
\par  \vskip 0.15 cm
(B) \ For each $ D_{i} \ ( 1\leq i \leq n ), $ we have \\
(1) \ $ C_{D_{i}}^{\prime} ( \Q _{2}) \neq \emptyset $ if and only
if one of the following conditions holds: \\
(a) \ $ D_{i} \equiv 1 (\text{mod} 8 ); $ \quad (b) \ $ ( 1 - p
\widehat{D_{i}} ) ( 1 - q \widehat{D_{i}} ) \equiv 0 (\text{mod} 16 ); $ \\
(c) \ $ D_{i} \equiv 3 (\text{mod} 8 ) $ and $ p \equiv 1
(\text{mod} 4 ); $ \quad (d) \ $ D_{i} \equiv 7 (\text{mod} 8 ) $
and $ p \equiv 3 (\text{mod} 4 ). $ \\
(2) \ $ C_{D_{i}}^{\prime} ( \Q _{p}) \neq \emptyset $ and $
C_{D_{i}}^{\prime} ( \Q _{q}) \neq \emptyset . $ \\
(3) \ For each $ j \neq i, \quad  C_{D_{i}}^{\prime} ( \Q
_{D_{j}}) \neq \emptyset \Longleftrightarrow \ \left( 1 -
\left(\frac{D_{i}}{D_{j}}\right) \right) \left( 1 - \left(\frac{p
q D_{i}}{D_{j}}\right) \right) = 0. $ \\
(4) \ $ C_{D_{i}}^{\prime} ( \Q _{D_{i}}) \neq \emptyset
\Longleftrightarrow \ \left( 1 - \left( \frac{ p
\widehat{D_{i}}}{D_{i}}\right) \right) \left( 1 - \left(\frac{q
\widehat{D_{i}}}{D_{i}}\right) \right) = 0. $
\par  \vskip 0.15 cm
(C) \ (1) \ $ C_{D}^{\prime } ( \Q _{2}) \neq \emptyset $ if and
only if one of the following conditions holds: \\
(a) \ $ D \equiv 1 (\text{mod} 8 ); $ \quad (b) \ $ p \equiv 1, 7
(\text{mod} 8 ); $ \\
(c) \ $ D \equiv 3 (\text{mod} 8 ) $ and $ p \equiv 5 (\text{mod}
8 ); $ \quad (d) \ $ D \equiv 7 (\text{mod} 8 ) $
and $ p \equiv 3 (\text{mod} 8 ). $ \\
(2) \ $ C_{D}^{\prime} ( \Q _{p}) \neq \emptyset $ and $
C_{D}^{\prime} ( \Q _{q}) \neq \emptyset . $ \\
(3) \ For each $ i, \quad  C_{D}^{\prime} ( \Q _{D_{i}}) \neq
\emptyset \Longleftrightarrow \ \left( 1 -
\left(\frac{p}{D_{i}}\right) \right) \left( 1 - \left(\frac{ q
}{D_{i}}\right) \right) = 0. $

\par  \vskip 0.2 cm

{\bf Proof. } \ The proof is similar to Proposition 2.3, we omit the
details.

\par  \vskip 0.2 cm

{\bf Corollary 2.8 } \ We assume $ \varepsilon = - 1 $ and the
elliptic curve $ E^{\prime} = E^{\prime}_{-} $ be as in (1.2).
\par  \vskip 0.15 cm
For each $ D_{i} \ ( 1\leq i \leq n ), $ denote $$ \Pi
 _{i}^{-}(D)^{\prime } = \left( 1 - \left(\frac{ p
\widehat{D_{i}}}{D_{i}}\right) \right) \left( 1 - \left(\frac{q
\widehat{D_{i}} }{D_{i}}\right) \right) + \sum_{j = 1, \ j \neq
i}^{n} \left( 1 - \left(\frac{D_{i}}{D_{j}}\right) \right) \left(
1 - \left(\frac{p q D_{i} }{D_{j}}\right) \right). $$ \ Then \ $
D_{i} \in S^{(\widehat{\varphi })} (E^{\prime } / \Q ) $ \ if
and only if one of the following conditions holds: \\
(1) \ $ \left \{
   \begin{array}{l}
  D_{i} \equiv 1 (\text{mod} 8 ), \\
  \Pi _{i}^{-}(D)^{\prime } = 0;
  \end{array}
  \right.   \qquad \quad (2) \ \left \{ \begin{array}{l}
  (1 - p \widehat{D_{i}} )(1 - q \widehat{D_{i}} )  \equiv 0 (\text{mod} 16 ), \\
  \Pi _{i}^{-}(D)^{\prime } = 0;
  \end{array}
  \right. $ \\
(3) \ $ \left \{
   \begin{array}{l}
  D_{i} \equiv 3 (\text{mod} 8 ), \\
  p \equiv 1 (\text{mod} 4 ), \\
  \Pi _{i}^{-}(D)^{\prime } = 0;
  \end{array}
  \right.   \qquad  \quad (4) \ \left \{ \begin{array}{l}
  D_{i} \equiv 7 (\text{mod} 8 ), \\
  p \equiv 3 (\text{mod} 4 ), \\
  \Pi _{i}^{-}(D)^{\prime } = 0.
  \end{array}
  \right. $
\par  \vskip 0.2 cm

{\bf Proof. } \ The results follow easily from Proposition 2.7.
\quad $ \Box$

Now we come to prove the Theorems.
\par  \vskip 0.2 cm

{\bf Proof of Theorem 1.1.} \ (1) \ For each $ i \in \{ 1, \cdots ,
n \}, $ \ we have $ \left[ \frac{1}{1 + \Pi _{i}^{+}(D)} \right] =
0, \ 1. $ \ Furthermore, \ $ \left[ \frac{1}{1 + \Pi _{i}^{+}(D)}
\right] = 1 \ \Longleftrightarrow \ \Pi _{i}^{+}(D) = 0 \
\Longleftrightarrow \ \left(\frac{-1}{D_{i}} \right) = 1 $ \ and \ $
\left(\frac{p}{D_{i}} \right) = \left(\frac{q}{D_{i}} \right) =
\left(\frac{D_{j}}{D_{i}} \right) = 1 \ ( 1 \leq j \leq n \
\text{and} \ j \neq i ) \ \Longleftrightarrow \ D_{i} \in
S^{(\varphi)} ( E / \Q ) $ \ by Proposition 2.1. Then the conclusion
(1) follows. \\
(2) \ Conclusion (2) follows easily from the above (1) and Corollary
2.2 (A). \ This proves Theorem 2.1. \quad $ \Box $

\par  \vskip 0.2 cm

{\bf Proof of Corollary 1.2.} \ By Proposition 2.1, we know that $ d
\notin S^{(\varphi )} (E / \Q ) $ for any $ d \in \Q (S, 2) $
satisfying one of the following conditions holds: (1) \ $ d < 0; $ \
(2) \ $ p \mid d; $ \ (3) \ $ q \mid d. $ Since $ S^{(\varphi )} (E
/ \Q ) \subset \Q ( S, 2 ) = < -1, 2, p, q, D_{1}, \cdots , D_{n} >,
$ we get $ S^{(\varphi )} (E / \Q ) \subset < 2, D_{1}, \cdots ,
D_{n} >. $ Furthermore, if condition (A) holds, then by Proposition
2.1, we have $ D_{1}, \cdots , D_{n} \in S^{(\varphi )} (E / \Q ), $
so $ < D_{1}, \cdots , D_{n} > \ \subset \ S^{(\varphi )} (E / \Q )
\ \subset \ < 2, D _{1}, \cdots , D_{n} >. $ If condition (B) holds,
then by Proposition 2.1, we have $ D_{1}, \cdots , D_{n} \in
S^{(\varphi )} (E / \Q ) $ but $ 2 \notin S^{(\varphi )} (E / \Q ),
$ hence $ S^{(\varphi )} (E / \Q ) \ \cong \ < D _{1}, \cdots ,
D_{n} > \cong \left( \Z / 2 \Z \right)^{ n }. $ If condition (C)
holds, then by Proposition 2.1, we have $ 2, D_{1}, \cdots , D_{n}
\in S^{(\varphi )} (E / \Q ), $ hence $ S^{(\varphi )} (E / \Q ) \
\cong \ < 2, D _{1}, \cdots , D_{n} > \cong \left( \Z / 2 \Z
\right)^{ n + 1 }. $ This proves Corollary 1.2. \quad $ \Box $
\par  \vskip 0.2 cm

{\bf Proof of Theorem 1.5.} \ The results of the Selmer groups are
given in Corollaries 1.2 and 1.4 above. \\
Now $ E(\Q )[2] = \{O, (0, 0), (-\varepsilon pD, 0), (-\varepsilon
qD, 0) \}, \ \varphi (E(\Q)[2]) = \{O, (0, 0) \} =
E^{'}(\Q)[\widehat{ \varphi }], $ by the exact sequences ([S] p.298,
314, 301 ) $$  \begin{array}{l} 0 \longrightarrow \frac{E^{'}(\Q
)}{\varphi (E(\Q ))} \longrightarrow S^{(\varphi )}(E /\Q )
\longrightarrow \text{TS}(E / \Q)[\varphi ]
\longrightarrow 0, \\
0 \longrightarrow \frac{E(\Q)}{\widehat{\varphi }(E^{\prime }(\Q ))}
\longrightarrow S^{( \widehat{\varphi })}(E^{\prime }/ \Q )
\longrightarrow \text{TS}(E^{\prime } / \Q)[\widehat{\varphi } ]
\longrightarrow 0, \\
0 \longrightarrow \frac{E^{'}(\Q)[\widehat{ \varphi }]}{\varphi
(E(\Q)[2])} \longrightarrow \frac{E^{'}(\Q)}{\varphi (E(\Q) ) }
\longrightarrow \frac{E(\Q)}{2E(\Q)} \longrightarrow
\frac{E(\Q)}{\widehat{ \varphi }(E^{'}(\Q))} \longrightarrow  0,
\quad \text{we get}
\end{array} $$
$$  \begin{array}{l}  \ \text{rank}(E(\Q)) + \dim_{2}(\text{TS}(E / \Q
)[\varphi ] ) + \dim_{2}(\text{TS}(E^{'} / \Q )
[\widehat{\varphi }]) \\
= \dim_{2}(S^{(\varphi )}(E / \Q ) ) + \dim_{2}(S^{(\widehat{\varphi
} )}(E^{'} / \Q ) ) - 2.
\end{array} $$
From this and Corollaries 1.2 and 1.4, the results in Theorem 1.5
follow, and the proof is completed. \quad $ \Box $
\par  \vskip 0.2 cm

Depending on the corresponding results in Propositions 2.3, 2.5,
2.7, the Theorems 1.3, 1.6, 1.8, 1.10 and their Corollaries can be
similarly proved.

\par  \vskip 0.4 cm

\hspace{-0.8cm} {\bf References }
\begin{description}

\item[[C]] J.W.S. Cassels,  Lectures on Elliptic Curves,  LMS 24,
Cambridge University Press, New York, 1995.

\item[[IK]] H. Iwaniec and E. Kowalski, Analytic Number Theory,
American Mathematical Society, Providence, Rhode Island, 2004.

\item[[KS]] R. Kloosterman and E. F. Schaefer, Selmer groups of
elliptic curves that can be arbitrarily large, Journal of Number
Theory, 99: 148-163 (2003).

\item[[L]] S. Lang,  Algebra, Revised Third Edition, GTM 211,
Springer-Verlag, New York, 2002.

\item[[MSS]] J.R. Merriman, S. Siksek and
N.P.Smart, Explicit 4-descents on an elliptic curve, Acta
Arithmetica lxxvii.(4):385-404(1996).

\item[[QZ]] D. Qiu, X. Zhang, Mordell-weil groups and selmer
groups of two types of elliptic curves, Science in China (series A),
2002, Vol.45, No.11, 1372-1380.

\item[[Sch]] S. Schmitt, Computation of the Selmer groups of certain
parametrized elliptic curve, Acta Arithmetica
lxxviii.(3):241-254(1997).

\item[[S]] J.H.Silverman,  The Arithmetic of Elliptic Curves,  GTM 106,
Springer-Verlag, New York, 1986.

\item[[ST]] R.J.Strocker and J.Top, On the equation $y^{2}=(x+p)(x^{2}+p^{2}), $
Rocky Mountain J. of Math. 24(1994), 1135-1161.

\end{description}

\par  \vskip 0.3 cm

\end{document}